\documentclass[12pt]{amsart}
\usepackage{epsf,overpic,amssymb,amscd,times,stmaryrd}

\newtheorem{thm}{Theorem}[section]

\newtheorem{lemma}[thm]{Lemma}

\newtheorem{q}[thm]{Question}

\numberwithin{equation}{subsection}

\newcommand{\R}{\mathbb{R}}
\newcommand{\Z}{\mathbb{Z}}

\newcommand{\F}{\mathcal{F}}

\newcommand{\bdry}{\partial}
\newcommand{\s}{\vskip.1in}
\newcommand{\n}{\noindent}

\newcommand{\be}{\begin{enumerate}}
\newcommand{\ee}{\end{enumerate}}

\begin{document}

\title{On the flux of pseudo-Anosov homeomorphisms}

\author{Vincent Colin}
\address{Universit\'e de Nantes, UMR 6629 du CNRS, 44322 Nantes, France}
\email{Vincent.Colin@math.univ-nantes.fr}

\author{Ko Honda}
\address{University of Southern California, Los Angeles, CA 90089}
\email{khonda@usc.edu} \urladdr{http://rcf.usc.edu/\char126 khonda}

\author{Fran\c cois Laudenbach}
\address{Universit\'e de Nantes, UMR 6629 du CNRS, 44322 Nantes, France}
\email{francois.laudenbach@univ-nantes.fr}

\date{This version: September 19, 2008.}



\thanks{VC supported by the Institut
Universitaire de France and the ANR Symplexe. KH supported by an NSF CAREER Award
(DMS-0237386) and NSF Grant DMS-0805352.}

\begin{abstract}
We exhibit a pseudo-Anosov homeomorphism of a surface $S$ which acts
trivially on $H_1 (S;\Z )$ and whose flux is nonzero.
\end{abstract}

\maketitle

\section{Introduction} \label{section: intro}

Let $S$ be a compact oriented surface with nonempty boundary,
$\omega$ be an area form on $S$, and $h$ be an area-preserving
diffeomorphism of $(S,\omega)$.  Consider the mapping torus $\Sigma
(S,h)$ of $(S,h)$, which we define as $(S\times [0,1])/(x,1) \sim
(h(x),0)$. Here $(x,t)$ are coordinates on $S\times[0,1]$. If there
is a contact form $\alpha$ on $\Sigma(S,h)$ for which $d\alpha
\vert_{S\times \{0\}} =\omega$ and the corresponding Reeb vector
field $R_\alpha$ is directed by $\bdry_t$, then we say $h$ is {\em
the first return map of $R_\alpha$}. In this note we investigate the
following question:

\begin{q} \label{question1}
What is the difference between an area-preserving diffeomorphism $h$
of a surface $(S,\omega)$ and the first return map of a Reeb flow
$R_\alpha$, defined on $\Sigma(S,\omega)$?
\end{q}

One easily computes that the first return map of $R_\alpha$ is
$\omega$-area-preserving (cf.\ Lemma~\ref{lemma: form-preserving}).
Question~\ref{question1} can then be rephrased as follows:

\begin{q} \label{question2}
Can every area-preserving $h$ be expressed as the first return map
of a Reeb flow $R_\alpha$?
\end{q}

We emphasize that we are interested in the rigid problem of
realizing a {\em given diffeomorphism $h$}, instead of its
realization {\em up to isotopy}. This question is of particular
importance when one tries to compute the {\it contact homology} of
a contact structure adapted to an open book decomposition
\cite{CH2}. The periodic orbits of an adapted Reeb flow that are
away from the binding of the open book correspond to periodic
points of the first return map. Hence we would like to understand
which monodromy maps can be realized by first return maps of Reeb
flows.

It turns out that the answer to Question~\ref{question2} is
negative. There is an invariant of an area-preserving
diffeomorphism $h$, called the {\em flux}, which is an obstruction
to $h$ being the first return map of a Reeb flow. In
Section~\ref{section:flux} we define the flux and also show that
it is easy to modify the flux of a diffeomorphism within its
isotopy class.

The case of particular interest to us is when $h$ is pseudo-Anosov.
Recall that a homeomorphism $h:S\stackrel\sim\rightarrow S$ is {\em
pseudo-Anosov} if there exist $\lambda>1$ and two transverse
singular measured foliations
--- the stable measured foliation $(\mathcal{F}^s,\mu^s)$ and the
unstable measured foliation $(\mathcal{F}^u, \mu^u)$ --- such that
$h(\mathcal{F}^s,\mu^s)=(\mathcal{F}^s,{1\over \lambda}\mu^s)$ and
$h(\mathcal{F}^u,\mu^u)=(\mathcal{F}^u,\lambda\mu^u)$. The
homeomorphism $h$ is a diffeomorphism away from the singular points
of the measured foliations. A pseudo-Anosov representative $h$ of a
mapping class is unique in the sense that any two pseudo-Anosov
homeomorphisms $h_1$, $h_2$ in the same mapping class are conjugate
via an everywhere smooth diffeomorphism $\phi$ which is isotopic to
the identity. In particular, such a $\phi$ sends the stable
foliation of $h_1$ to the stable foliation of $h_2$ and the unstable
foliation of $h_1$ to the unstable foliation of $h_2$.
(See~\cite[Expos\'e 12, Th\'eor\`eme III and Lemma 16 for
smoothness]{FLP}.) We define the area form $\omega$ to be given by
the product of $\mu^s$ and $\mu^u$. The form $\omega$ is the unique
$h$-invariant area form up to a constant multiple, and is singular
in the sense that it vanishes at the singular points of the
invariant foliations. Now, the pseudo-Anosov case is of special
interest since the pseudo-Anosov homeomorphism is a rigid
representative in its mapping class (hence the flux can be seen as
an invariant of the mapping class) and also since it is known that
every contact structure is carried by an open book decomposition
whose monodromy is isotopic to a pseudo-Anosov
homeomorphism~\cite{CH1}. Hence we ask the following question:

\begin{q}
Can every pseudo-Anosov homeomorphism $h$ be expressed as the first
return map of a Reeb flow $R_\alpha?$
\end{q}

The main theorem of this paper is Theorem~\ref{thm:flux}, which
states that the answer to this question is also negative, i.e., the
flux is not always zero for pseudo-Anosov homeomorphisms.

\section{The flux}
\label{section:flux}

The goal of this section is to give basic properties of the flux;
see \cite{Ca}. The discussion will be done more generally on a
compact symplectic manifold, since it might be more transparent in
that context.

\subsection{Flux}

Let $(S,\omega)$ be a compact symplectic manifold and $h$ be a
symplectomorphism of $(S,\omega)$. Let $h_* :H_1 (S;\Z )\rightarrow
H_1 (S;\Z )$ be the map on homology induced from $h$ and let $K$ be
the kernel of $h_* -id$.  Also let $\Gamma$ be a lattice of $\R$
generated by $\int_\Sigma \omega$, where $[\Sigma]$ ranges over
$H_2(S;\Z)$.  Then define the map
$$F_h:K\rightarrow \R/\Gamma$$
as follows: Let $[\gamma]\in K$. Since $\gamma$ is homologous to
$h(\gamma)$, one can find an oriented singular cobordism $C$
(mapped into $S$) whose boundary consists of $h(\gamma)-\gamma$.
We then define
$$F_h(\gamma)=\int_C\omega.$$
Two cobordisms $C,C'$ with the same boundary differ by an element of
$H_2(S;\Z)$; hence the quantity is well-defined only up to $\Gamma$.
It is straightforward to verify that $F_h(\gamma)$ also only depends
on the homology class of $\gamma$. The number $F_h ([\gamma])\in
\R/\Gamma$ is thus well-defined and is called the {\it flux of $h$
along $\gamma$}.  We say the flux of $h$ is {\em nonzero} if the
image of $K$ is not $[0]\in \R/\Gamma$.

If $h_1,h_2$ are two symplectomorphisms of $(S,\omega)$ and
$[\gamma]=[h_1(\gamma)]=[h_2(\gamma)]$, then $$F_{h_2}\circ
F_{h_1}([\gamma])=F_{h_2}([\gamma])+F_{h_1}([\gamma]).$$ In other
words, the flux is a homomorphism, when viewed as a map from the
group $Symp_0(S,\omega)$ of symplectomorphisms which act trivially
on $H_1(S;\Z)$ to $Hom(H_1(S;\Z),\R/\Gamma)=H^1(S;\R/\Gamma).$  We
can also easily modify the flux of any $h\in Symp_0(S,\omega)$ by
composing with time-1 maps of locally Hamiltonian flows.

If in addition $\omega=d\beta$,
then the form $h^* \beta -\beta$ is a
closed $1$-form and the flux of $h$ along $\gamma$ can be
rewritten as
$$F_h([\gamma])=\int_{\gamma} h^*\beta -\beta,$$
by the use of Stoke's formula. The flux of $h$ is nonzero if and
only if $[h^* \beta -\beta ]\neq 0$ on $K$.  Moreover, $\Gamma=0$.

\subsection{$2$-forms on the mapping torus}
Let $\Sigma(S,h)=(S\times[0,1])/(x,1)\sim (h(x),0)$ be the mapping
torus of $(S,\omega)$. It fibers over the circle with fiber $S$.

There is a natural closed $2$-form $\omega_h$ on $\Sigma (S,h)$,
which is obtained by setting $\omega_h=\omega$ on $S\times[0,1]$
and identifying via the symplectomorphism $h$. The $2$-form
$\omega_h$ pulls back to $\omega$ on $S\times \{ t\}$, $t\in
[0,1]$, and its kernel is directed by $\partial_t$, where $t$ is
the coordinate for $[0,1]$.

We have the following lemmas:

\begin{lemma} \label{lemma: form-preserving}
Suppose $\omega$ is exact. If $h$ is the first return map of a
Reeb vector field $R_\alpha$ where $\alpha$ satisfies $d\alpha
\vert_{S\times \{ 0\}} =\omega$, then $h$ is a symplectomorphism
of $(S,\omega)$. Moreover, $d\alpha=\omega_h$.
\end{lemma}

\begin{proof}
Consider the contact $1$-form $\alpha=fdt+\beta$ on $S\times[0,1]$,
where $f=f(x,t)$ is a function and $\beta=\beta(x,t)$ is a $1$-form
in the $S$-direction.  Write $\beta_t(x)=\beta(x,t)$.  We compute
$$d\alpha=d_Sf\wedge dt+d_S\beta_t+dt\wedge\dot\beta_t,$$ where
$d_S$ is the exterior derivative in the $S$-direction and
$\dot\beta_t={d\beta_t\over dt}$.  By the condition
$d\alpha|_{S\times\{0\}}=\omega$, we have $d_S\beta_0=\omega$. Since
we can normalize $R_\alpha=g\bdry_t$, where $g=g(x,t)$, it follows
that $i_{R_\alpha}d\alpha= g(-d_Sf+\dot\beta_t)=0$ and
$\dot\beta_t=d_Sf$ is an exact form on $S$. Hence $d_S\beta_t$ is
independent of $t$ and equals $\omega$. This shows that
$d\alpha=d_S\beta_t=\omega_h$.  By the invariance of $\alpha$ under
the map $(x,t)\mapsto (h(x),t-1)$, we see that $h$ preserves
$\omega$.
\end{proof}

\begin{lemma}\label{lemma:monotonicity}
Suppose $\omega$ is exact and the flux of $h$ is nonzero.  Then
$[\omega_h]$ is nonzero in $H^2(\Sigma(S,h);\Z)$.  Hence $h$ cannot
be realized as the first return map of a Reeb vector field
$R_\alpha$.
\end{lemma}

\begin{proof}
Let $\gamma$ be a curve in $S$ such that $F_h (\gamma)$ is
nonzero. Then $-\gamma$ and $h(\gamma )$ bound a subsurface
$C\subset S\times\{0\}$ so that $\int_C \omega \neq 0$. We
construct a closed $2$-cycle $C'$ in $\Sigma (S,h )$ by gluing
$\gamma \times [0,1]$ with $C$. Now we see that $\int_{C'}
\omega_h =\int_C \omega \neq 0$. Hence $[\omega_h]\not=0$ in
$H^2(\Sigma(S,h);\Z)$. By Lemma~\ref{lemma: form-preserving}, $h$
cannot be the first return map of a Reeb vector field.
\end{proof}

Conversely, when $\omega$ is exact and $h$ is the identity near
$\partial S$, we have a criterion, due to Giroux (see~\cite{Co}), to
realize $h$ as the first return map of a Reeb vector field. The
condition that $h=id$ near $\bdry S$ is not realized in general for
pseudo-Anosov homeomorphisms, but in practice it is possible to
deform the diffeomorphism near $\partial S$ so that it is the
identity, without altering the sets of periodic points too much;
see~\cite{CH2}.

\begin{lemma}[Giroux]\label{lemma: Reeb constr}
Let $(S,\omega=d\beta)$ be a compact exact symplectic manifold and
$h$ be a symplectomorphism of $(S,\omega)$, which is the identity
near $\partial S$. If $[h^* \beta -\beta ]=0$ in $H^1(S;\R)$, then
there exists a contact form $\alpha$ on $\Sigma (S,h)$ and a Reeb
vector field $R_\alpha$ whose first return map on one fiber is $h$.
\end{lemma}

\begin{proof}
We have that $h^* \beta -\beta =df$. Note that $df=0$ near $\bdry S$
since $h=id$ near $\bdry S$. One can always translate $f$ so that
$f$ is strictly positive on $S$ and is constant near $\partial S$.
The $1$-form $\alpha =dt+\beta$ is a contact form on $S\times \R$
whose Reeb vector field is $\partial_t$. It is invariant under the
diffeomorphism
$$H : (x,t)\mapsto (h (x),t-f(x)),$$
and thus induces a contact form $\alpha$ on $\Sigma (S,h)\simeq
(S\times \R )/((x,t)\sim H (x,t))$.
\end{proof}

\subsection{Surface case}
Let us now specialize to the case of interest: $S$ is a compact
oriented surface, $\omega$ is an area form on $S$, and $h$ is an
area-preserving diffeomorphism of $(S,\omega)$. Let us assume
without loss of generality that the $\omega$-area of $S$ is $1$. If
$\bdry S=\emptyset$, then $\Gamma=\Z$ since $H_2(S;\Z)$ is generated
by $[S]$.  On the other hand, if $\bdry S\not=\emptyset$, then
$\Gamma=0$ and the exactness of $\omega$ is automatically satisfied.

The goal of this note is to prove the following:

\begin{thm}\label{thm:flux}
There exist a compact surface $S$ with empty (resp. non-empty)
boundary and a pseudo-Anosov homeomorphism $h$ of $S$ with $h_*
=id$, whose flux with respect to the singular $h$-invariant area
form $\omega$ is nonzero, as viewed in $\R/\Z$ (resp.\ $\R$).
\end{thm}

We now discuss a technical issue, namely the fact that $h$ is only
$C^0$ at the singular set $L=\{p_1,\dots,p_k\}$ of the
stable/unstable foliations. Let $\omega $ be the $h$-invariant
singular area form given by the product $\mu^u\otimes \mu^s$ of both
transverse measures. It is singular in the sense that it is a
$2$-form which vanishes on $L$. As a measure it is equivalent to any
Lebesgue measure on $S$. Hence, according to a theorem of
Oxtoby-Ulam \cite{OU}, it is conjugated by a homeomorphism to a
smooth area form.

Instead of the Oxtoby-Ulam approach, our approach will be based on
Moser's lemma. Let $D$ be an arbitrarily small open neighborhood of
$L$ so that each connected component of $D$ is a polygonal region
whose boundary consists of subarcs of leaves of $\mathcal F^s$ or
$\mathcal F^u$.  Then we have the following:

\begin{lemma}
There exist an everywhere smooth area form $\omega'$ on $S$ and a
diffeomorphism $h'$, which coincide respectively with $\omega$ and
$h$ outside of $D$ and satisfy (i)
$\int_{D_0}\omega=\int_{D_0}\omega'$ for each connected component
$D_0$ of $D$ and (ii) $(h')^*\omega'=\omega'$.
\end{lemma}

\begin{proof}
Let $\omega'$ be an area form which coincides with $\omega$ on
$S-(D\cap h(D))$, and has the same area as $\omega$ on each
connected component $D_0$ of $D$. (By using an auxiliary area form
on $S$, the construction of such an $\omega'$ becomes equivalent to
the extension of a positive smooth function with a fixed integral.)
There also exists a smooth diffeomorphism $\psi$ of $S$ which
coincides with $h$ on $S-D$. Note that the germ of $h$ along
$\partial D$ extends to an embedding of $D$ into $S$. By the
construction of $\omega'$ and $\psi$, we have
$\psi^*\omega'=\omega'$ on $S-D$.

We now claim that
$$ \int_{D_0}\psi^*\omega'= \int_{D_0}\omega'$$
for each component $D_0$ of $D$. We have
$\int_{D_0}\omega'=\int_{D_0}\omega=\int_{D_0}h^*\omega=\int_{h(D_0)}\omega,$
by our choice of $\omega'$ and the $h$-invariance of $\omega$.  On
the other hand, we have $\int_{D_0}\psi^*\omega'=
\int_{h(D_0)}\omega'$ by a change of variables. If $D_0'$ is the
component of $D$ that nontrivially intersects $h(D_0)$, then
$\int_{D_0'\cap h(D_0)}\omega=\int_{D_0'\cap h(D_0)}\omega'$, since
$\int_{D_0'}\omega=\int_{D_0'}\omega'$ and $\omega=\omega'$ on
$D_0'-h(D_0)$. From this we deduce that
$\int_{h(D_0)}\omega=\int_{h(D_0)}\omega'$. The claimed equality
follows.

Finally, Moser's lemma applies on $D$ to the pair of area forms
$\omega'$ and $\psi^*\omega'$. It yields a diffeomorphism $\varphi$
of $D$ which is the identity near the boundary (hence extends to $S$
by the identity of $S'$) such that
$\varphi^*(\psi^*\omega')=\omega'$. We set $h'=\psi\circ \varphi$.
This diffeomorphism meets the required condition both on $S'$ and
$D$, hence on $S$.
\end{proof}

If we choose $\gamma$ so that both $\gamma$ and $h'(\gamma)$ avoid
the small neighborhood $D$ of the singular locus $L$ (after
isotopy), then we see that $F_h(\gamma)=F_{h'}(\gamma)$. Since the
flux only depends on the curve up to isotopy, it follows that
$F_h=F_{h'}$.

\s\n {\bf Remark.} When $\F^s$ and $\F^u$ are orientable, the
transverse measures define $1$-forms that are closed but not exact.
They are eigenvectors for $h^*$ with eigenvalues $\lambda$ and
${1\over\lambda}$. Thus if $h_*=id$, then the foliations are not
orientable.

\section{Proof of Theorem~\ref{thm:flux}}

Let $S=S_g$ be a closed oriented surface of genus $g$ and $\alpha$
and $\beta$ be two $1$-dimensional submanifolds of $S$, i.e., the
union of disjoint simple closed curves.

We recall that $\alpha$ and $\beta$ {\em fill $S$} if $\alpha$ and
$\beta$ intersect transversely and minimally and if each region of
$S- (\alpha \cup \beta )$ is a $2n$-gon with $n>1$. Such a system
of curves allows one to define two systems of flat charts, the
$\alpha$- and the $\beta$-charts, in the following way: The set
$\alpha \cup \beta$ gives a cell decomposition of $S$. Consider
its dual cell decomposition.  (By this we mean we place a vertex
$v_{P_i}$ in the interior of each component $P_i$ of $S-
(\alpha\cup\beta)$.  If $P_i$ and $P_j$ share an edge of
$\alpha\cup\beta$, then take an edge from $v_{P_i}$ to $v_{P_j}$
which passes through the common edge of $\alpha\cup\beta$ exactly
once.) Let $E_\beta$ be the union of edges of the dual cellular
decomposition that meet $\beta$. Then $E_\beta$ cuts $S$ into
annuli whose cores are the components of $\alpha$ that we call the
$\alpha$-charts.
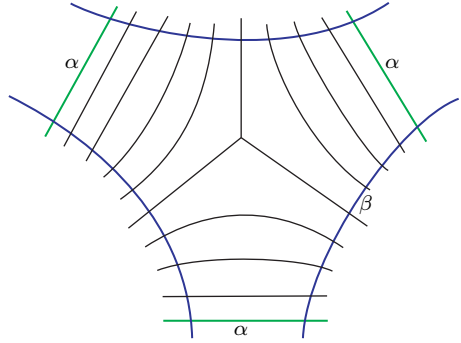
\begin{figure}[ht]
\begin{overpic}[width=6cm]{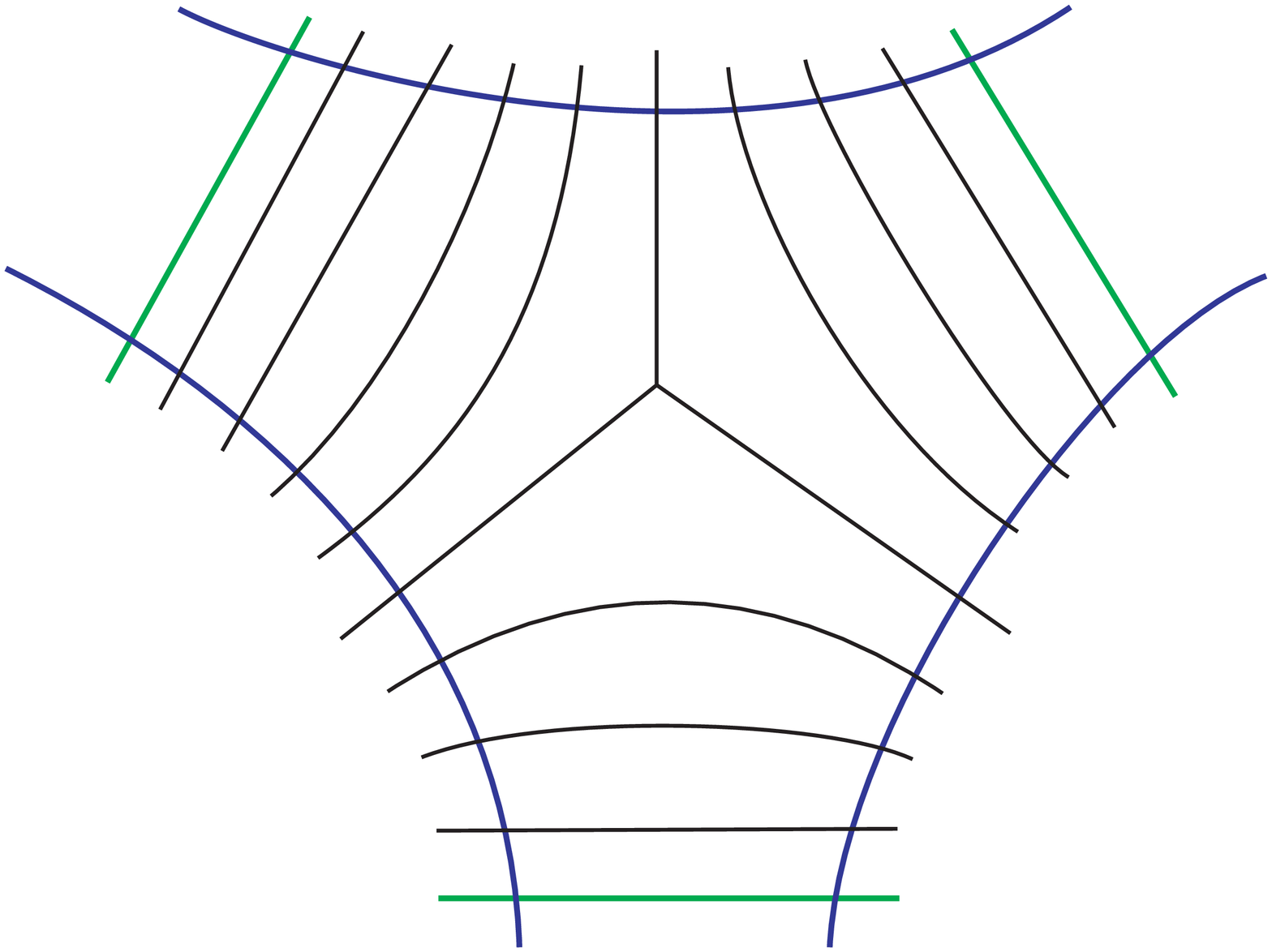}
\put(12.5,60){\tiny $\alpha$} \put(77.7,28.6){\tiny $\beta$}
\put(83.5,60){\tiny $\alpha$} \put(50,0.8){\tiny $\alpha$}
\end{overpic}
\caption{Dual cell decomposition meeting the $\beta$ curves,
together with flat geodesics parallel to $\alpha$.} \label{hexagon}
\end{figure}
The $\beta$-charts are defined similarly. Note that there is one
chart for each curve and hence each chart can be viewed as a
thickening of an appropriate $\alpha$- or $\beta$-curve. These
charts are equipped with a singular flat metric $g$ which is
standard on each little square, corresponding to intersections of
$\alpha$- and $\beta$-charts, as explained in \cite[Expos\'e 13,
Section III]{FLP}. (In particular, the $\alpha$-metric and the
$\beta$-metric coincide on the squares.)

We will construct our example on a surface $S_5$ of genus $5$.

\begin{lemma}\label{lemma:construction}
There exist two multicurves $\alpha =\alpha_1 \cup \alpha_2$ and
$\beta =\beta_1 \cup \beta_2$ on $S_5$ where:
\begin{itemize}
\item $\alpha$ and $\beta$ fill $S_5$; \item $\alpha_1$ and
$\beta_1$ are disjoint and form a bounding pair; \item $\alpha_2$
and $\beta_2$ are separating curves; \item $\# (\alpha_1 \cap
\beta_2 )=\# (\alpha_2 \cap \beta_1)=2$; \item $\# (\alpha_2 \cap
\beta_2 )=16$.
\end{itemize}
\end{lemma}

\begin{proof}
We start with a genus $2$ surface $S_2'$, together with simple
closed curves $\alpha_2'$ and $\beta_2'$ which are both
nullhomologous in $S_2'$, fill $S_2'$ and intersect $8$ times. See
Figure~\ref{genus2}.  (To see that $\beta_2'$ separates, take the
algebraic intersection number with a suitable basis for
$H_2(S_2';\Z)$.)
\begin{figure}[ht]
\begin{overpic}[width=9cm]{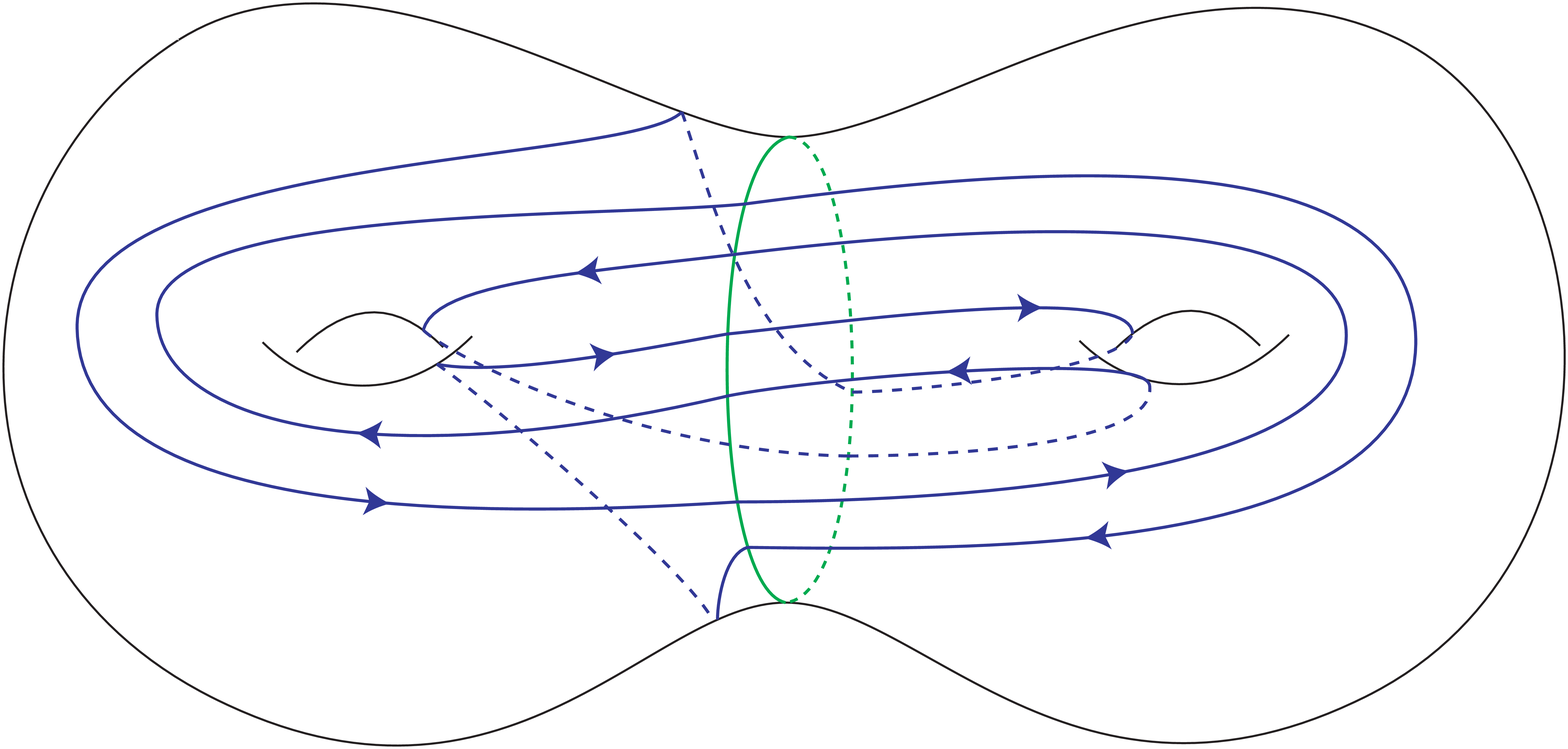}
\put(48.7,40.75){\small $\alpha_2'$} \put(70,38.25){\small
$\beta_2'$}
\end{overpic}
\caption{The genus $2$ surface $S_2'$.} \label{genus2}
\end{figure}
Now, two regions $H_1'$ and $H_2'$ of $S_2' \setminus (\alpha_2'
\cup \beta_2' )$ are $8$-gons.  For $i=1,2$, pick a disk $D_i'
\subset Int (H_i' )$.

We now take a second copy $(S_2'' ,\alpha_2'' ,\beta_2'' ,D_1'',
D_2'')$ of $(S_2' ,\alpha_2' ,\beta_2' ,D_1', D_2')$ and glue $S_2'
\setminus (D_1' \cup D_2' )$ to $S_2'' \setminus (D_1'' \cup D_2'')$
by identifying $\partial D_i'$ and $\partial D_i''$, $i=1,2$. We
call $S_5$ the resulting surface.  See Figure~\ref{genus5}.  Let
$\alpha_1 =\partial D_1' =\partial D_1''$ and $\beta_1 =\partial
D_2' =\partial D_2''$.

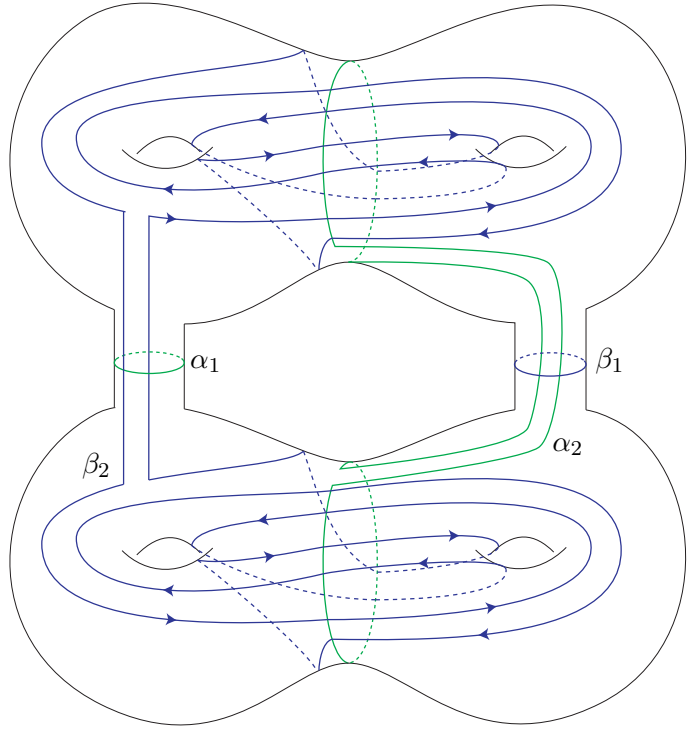
\begin{figure}[ht]
\begin{overpic}[width=9cm]{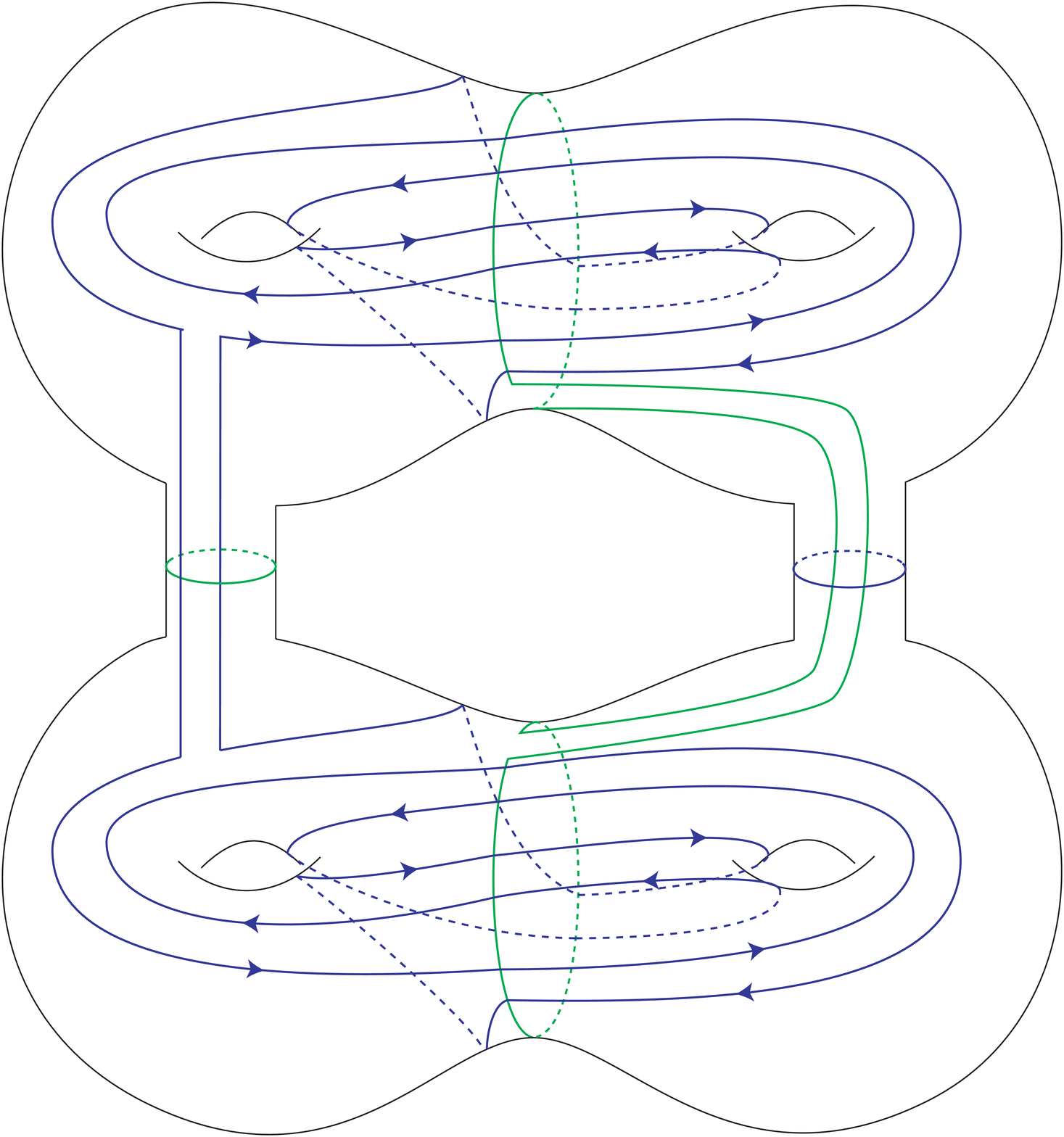}
\put(25,49.5){\small $\alpha_1$} \put(81,49.5){\small $\beta_1$}
\put(10,35){\small $\beta_2$} \put(75,38){\small $\alpha_2$}
\end{overpic}
\caption{The genus $5$ surface $S_5$.} \label{genus5}
\end{figure}

Next, take one connected component of $\beta_2' \cap \partial H_1'$
and one connected component of $\beta_2'' \cap \partial H_1''$, and
make the connected sum of these two components along an arc which
crosses $\alpha_1$ exactly once and stays inside $(H_1' \setminus
Int (D_1'))\cup (H_1'' \setminus Int (D_1'' ))$. We call $\beta_2$
the result of this sum of $\beta_2'$ and $\beta_2''$. By
construction, $\# (\alpha_1 \cap \beta_2 )=2$.  Now do the same
operation with components of $\alpha_2'$ and $\alpha_2''$ in
$\partial H_2'$ and $\partial H_2''$, so that the resulting curve
$\alpha_2$ satisfies $\# (\alpha_2 \cap \beta_1 )=2$.

By construction, we see that $\# (\alpha_2 \cap \beta_2 )=8+8=16$.
The families $\alpha =\alpha_1 \cup \alpha_2$ and $\beta =\beta_1
\cup \beta_2$ fill $S_5$. Since $\alpha_2'$ and $\alpha_2''$ were
nullhomologous, the same also holds for $\alpha_2$. Finally it is
clear that $\alpha_1$ and $\beta_1$ are disjoint and cobordant in
$S_5$.
\end{proof}

The system $\alpha=\alpha_1\cup\alpha_2$, provided by
Lemma~\ref{lemma:construction}, comes with two (oriented)
$\alpha$-charts $U_1 \supset \alpha_1$ and $U_2 \supset \alpha_2$,
where $U_i= [0,n_i ] \times [0,1]/(0,y)\sim (n_i ,y)$, $n_i$ denotes
$\# (\alpha_i \cap \beta)$, namely $n_1 =2$, $n_2 =18$, and
$\alpha_i=[0,n_i]\times\{{1\over 2}\}/\sim$. Similarly, there are
two $\beta$-charts $V_1$ and $V_2$, of the form $[0,1]\times
[0,m_i]/(x,0)\sim (x,m_i )$, where $m_i =\# (\beta_i \cap \alpha )$,
i.e., $m_1 =2$, $m_2 =18$,  and $\beta_i =\{ \frac{1}{2} \} \times
[0,m_i ]/\sim$. In what follows, we equip $S_5$ with the flat metric
associated to the system $\alpha$ and $\beta$ and compute areas
using this metric, normalized so that the total area of $S$ is $1$.

We will denote $[\tau_\eta]$ the mapping class of a positive Dehn
twist about the closed curve $\eta$. The class $[\tau_{\alpha_i}]$
admits an affine representative $\tau_{\alpha_i}$ which is given
on $U_i$ by the matrix
$$
\begin{pmatrix}
1&n_i\\
0&1
\end{pmatrix},$$
and is the identity on $U_j$ for $j\neq i$. Similarly,
$[\tau^{-1}_{\beta_i}]$ admits an affine representative
$\tau_{\beta_i}^{-1}$ which is given on $V_i$ by the matrix
$$
\begin{pmatrix}
1&0\\
m_i&1
\end{pmatrix},$$ and is the identity on $V_j$ for $j\neq i$.

\begin{lemma}
The map $h= \tau_{\alpha_2}  \circ \tau_{\alpha_1}^{9} \circ
\tau_{\beta_1}^{-9} \circ \tau_{\beta_2}^{-1}$ is a pseudo-Anosov
homeomorphism which acts by the identity on $H_1(S;\Z)$.
\end{lemma}

\begin{proof}
On both $U_1$ and $U_2$, the composition $\tau_{\alpha_2}  \circ
\tau_{\alpha_1}^9$ is given by the matrix
$$
\begin{pmatrix}
1&18\\
0&1
\end{pmatrix}$$
and thus is a smooth representative of its mapping class outside the
singular points of the flat structure. Similarly, on both $V_1$ and
$V_2$, the composition $\tau_{\beta_1}^{-9} \circ
\tau_{\beta_2}^{-1}$ is given by the matrix
$$
\begin{pmatrix}
1&0\\
18&1
\end{pmatrix}.$$
As a result, the homeomorphism $h$ is given away from the singular
points of the flat structure by the matrix
$$
\begin{pmatrix}
325&18\\
18&1
\end{pmatrix}. $$
Since the trace of the matrix is $>2$, $h$ is pseudo-Anosov. It
preserves the area coming from the singular flat metric on the
charts.

Since $\alpha_2$ and $\beta_2$ are homologous to zero and $\alpha_1$
and $\beta_1$ form a bounding pair, $h$ induces the identity on
homology.
\end{proof}

\begin{lemma} \label{lemma: nonzero flux}
The flux of $h$ is nonzero, when viewed in $\R/\Z$.  More precisely,
if $\gamma$ is a curve so that $\tau_{\beta_2}^{-1}(\gamma)$ has
geometric intersection one with each of $\alpha_1$ and $\beta_1$,
then $F_h ([\gamma])\neq 0$.
\end{lemma}

\begin{proof}
Let $\delta$ be a closed geodesic with respect to the singular flat
metric which corresponds to the singular flat coordinate system. In
other words, $\delta$ is a piecewise affine curve, with corners at
singularities of the affine structure.

First we claim that $\tau_{\beta_2}^{-1}$ has zero flux, i.e., the
area between $\delta$ and $\tau_{\beta_2}^{-1}(\delta)$ is zero for
all $\delta$. (Note that $\tau_{\beta_2}^{-1}(\delta)$ is not
necessarily a flat geodesic even if $\delta$ is.) Indeed, since
$\tau_{\beta_2}^{-1}$ is the identity on $V_1$, we only have to look
on $V_2$. The curve $\delta$ intersects $V_2$ along a finite union
of affine arcs $a_1 ,\dots ,a_n$. For any such $a_i$, the
concatenation $\overline{a}_i$ of $-a_i$ and $\tau_{\beta_2}^{-1}
(a_i)$ divides $V_2$ into two components with the same area. This
means that the area between $\overline{a}_i$ and $\beta_2 =\{
\frac{1}{2} \} \times [0,18]/\sim$ is zero. Thus, $\overline{a}_i$
bounds a subsurface in $S$ with the same area as the surface bounded
by $\beta_2 =\{ \frac{1}{2} \} \times [0,18]/\sim$. The sign of this
area depends on the sign of the intersection of $a_i$ with
$\beta_2$. Now observe that $\tau_{\beta_2}^{-1} (\delta)-\delta=
\cup_{1\leq i \leq n} \overline{a}_i$. Since $\beta_2$ is homologous
to zero, it has as many positive intersections with $\delta$ as
negative intersections. Thus the total signed area between $\delta$
and $\tau_{\beta_2}^{-1} (\delta)$ is the total signed area bounded
by $\cup_{1\leq i \leq n} \overline{a}_i$, which in turn is zero.
Similarly, we see that $\tau_{\alpha_2}$ has zero flux.

Next suppose the geodesic $\delta$ has geometric intersection one
with each of $\alpha_1$ and $\beta_1$.  We claim that the area
between $\delta$ and $\tau_{\alpha_1} \circ
\tau_{\beta_1}^{-1}(\delta)$ equals the area $A$ bounded by
$\alpha_1 \cup \beta_1$ in $S$. Since $\alpha_1 \cap \beta_1
=\emptyset$, it follows that $Int (U_1 )\cap Int (V_1)=\emptyset$
and the affine representatives $\tau_{\alpha_1}$ and
$\tau_{\beta_1}^{-1}$ commute. The curve $\delta$ intersects $V_1$
along a connected affine arc $b$ and $U_1$ along a connected affine
arc $a$. The concatenation of $-b$ and $\tau_{\beta_1}^{-1} (b)$ is
a closed curve $\overline{b}$ which cuts $V_1$ into two components
of the same area. Similarly, $\overline{a}$, obtained as the
concatenation of $-a$ and $\tau_{\alpha_1} (a)$, divides $U_1$ into
two components of the same area. Then $\tau_{\alpha_1} \circ
\tau_{\beta_1}^{-1} (\delta)-\delta$ equals $\overline{a} \cup
\overline{b}$, and $\delta$ and $\tau_{\alpha_1} \circ
\tau_{\beta_1}^{-1}(\delta)$ cobound a subsurface in $S_5$ of area
$A$. By the commutativity of $\tau_{\alpha_1}$ and
$\tau^{-1}_{\beta_1}$, we have $\tau_{\alpha_1}^9 \circ
\tau_{\beta_1}^{-9} =(\tau_{\alpha_1} \circ \tau_{\beta_1}^{-1}
)^9$.  Hence the area between $\delta$ and $\tau_{\alpha_1}^9\circ
\tau_{\beta_1}^{-9}(\delta)$ is $9A$.

We now claim that $A={1\over 2}\mbox{Area}(S)=\frac{1}{2}$. This is
due to the symmetry of the $\alpha$- and $\beta$-charts: The chart
$U_1$ is decomposed by $\alpha_1 = [0,2] \times
\{\frac{1}{2}\}/\sim$ into two pieces with the same area.  On the
other hand, $U_2$ is decomposed by $\beta_1$ into two rectangles
$R_1$ and $R_2$.  On each $R_i$, the number of intersections between
$\alpha_2$ and $\beta_2$ is $8$.  Hence
$\mbox{Area}(R_1)=\mbox{Area}(R_2)$. We conclude that the area
between $\delta$ and $\tau_{\alpha_1}^9 \circ
\tau_{\beta_1}^{-9}(\delta)$ is $9A=\frac{9}{2} \equiv
\frac{1}{2}\in \R/\Z$.

By putting together the above calculations and observing that
$F_h([\gamma])$ only depends on the isotopy class of $\gamma$, we
see that $F_h([\gamma])= A={1\over 2}$.
\end{proof}

This completes the proof of Theorem~\ref{thm:flux} when $\partial
S=\emptyset$.

To treat the case with boundary, we notice that the homeomorphism we
have constructed fixes the singular points of the invariant
foliations. We pick one of them and blow up the surface at this
point. The homeomorphism $h$ lifts to a pseudo-Anosov homeomorphism
on the blown-up surface $\overline{S}$ which fixes the blown-up
foliations. It also induces $id_*$ on $H_1 (\overline{S};\Z)$ and
has nonzero flux.

\s On the other hand, it is easy to construct pseudo-Anosov
homeomorphisms with vanishing flux that act trivially on $H_1
(S;\Z)$.  Let $\alpha_1$, $\beta_1$ be simple closed curves which
fill $S$ and are both nullhomologous. As explained in \cite[Expos\'e
13, Section III]{FLP}, if we compose twists along these curves
(positive Dehn twists along $\alpha_1$ represented by
$\tau_{\alpha_1}$ and negative Dehn twists along $\beta_1$
represented by $\tau^{-1}_{\beta_1}$, where we use at least one
$\tau_{\alpha_1}$ and at least one $\tau^{-1}_{\beta_1}$), we obtain
a pseudo-Anosov homeomorphism. The argument developed in the proof
of Lemma~\ref{lemma: nonzero flux} tells us that the flux is always
zero. More precisely, consider the singular flat metric compatible
with $\alpha_1$ and $\beta_1$, and let $\delta$ be a closed curve
represented by a flat geodesic.  As in the second paragraph of the
proof of Lemma~\ref{lemma: nonzero flux}, the area between $\delta$
and $\tau_{\alpha_1}(\delta)$ is zero for all $\delta$, since
$\alpha_1$ is separating.  Similarly, the area between $\delta$ and
$\tau_{\beta_1}^{-1}(\delta)$ is zero for all $\delta$. Hence, if
$h$ is any composition of $\tau_{\alpha_1}$ and
$\tau_{\beta_1}^{-1}$ (both with zero flux), then it also has zero
flux. This justifies the fact that, in the proof of
Theorem~\ref{thm:flux}, we have to look at more elaborate examples
to find nonzero flux.

Another case when $h$ can be realized as the first return map of a
Reeb flow is when $h^* -id$ is invertible and $\partial S\neq
\emptyset$.  We learned the following lemma from Yasha Eliashberg.

\begin{lemma}\label{lemma:surjective}
Let $h$ be a diffeomorphism of a surface $S$ with nonempty boundary
which preserves an area form $\omega$. If $1$ is not an eigenvalue
of $h^*$, then $\omega$ admits a primitive $\beta$ such that $[h^*
\beta -\beta ] =0$ in $H^1(S;\R)$.
\end{lemma}

\begin{proof}
Pick any primitive $\beta_0$ of $\omega$. By hypothesis, the map
$h^*-id$ is surjective. Thus, one can find $[\theta ]\in H^1(S;\R)$
such that $[h^* \beta_0 -\beta_0 ]=(h^*-id)[\theta]$. Now we have
that $\beta =\beta_0 -\theta$ is a primitive of $\omega$ and that
$[h^* \beta -\beta ] =0$ in $H^1(S;\R)$.
\end{proof}

Now, by applying Lemma~\ref{lemma: Reeb constr}, $h$ can be
realized as the first return map of a Reeb vector field.

\s
We end this section with the following questions:

\begin{q}
Is it possible to find a pseudo-Anosov homeomorphism of a surface
$S$ which acts trivially on $H_1 (S;\Z)$ and takes some
noncontractible curve $\gamma$ to a curve $h(\gamma )$ that can be
isotoped away from $\gamma$?\footnote{Dan Margalit has informed us
of an example of a pseudo-Anosov homeomorphism on a genus $3$
surface with this property.  His example also would therefore also
have nonzero flux.}
\end{q}

If yes, the flux of such a pseudo-Anosov homeomorphism would
automatically be nonzero.

\begin{q}
Let $g$ and $h$ be two pseudo-Anosov homeomorphisms acting trivially
on $H_1 (S;\Z )$ such that the composition $g\circ h$ is isotopic to
a pseudo-Anosov homeomorphism $f$. Suppose the flux of $g$ is zero
and the flux of $h$ is nonzero. Is the flux of $f$ nonzero?
\end{q}

If yes, this procedure would allow us to produce many pseudo-Anosov
homeomorphisms with nonzero flux.

\section{A question}

There is an invariant of an isotopy class of surface diffeomorphisms
$[h]$ which is defined in a manner much like the flux. We thank Ian
Agol for bringing this to the authors' attention. Let $S$ be a
hyperbolic surface with geodesic boundary. If $[\gamma]\in K$, i.e.,
$[h(\gamma)-\gamma]=0$, then represent $h(\gamma)$ and $\gamma$ by
geodesics, and compute the area bounded by the two geodesics.  By
the Gauss-Bonnet theorem, this area equals $-2\pi \chi(A)$, where
$A$ is a surface between the two geodesics. Here the Euler
characteristic $\chi(A)$ is more precisely an {\em Euler measure},
i.e., it is computed with signs: if $-A$ denotes $A$ with reversed
orientation, then one has $\chi (-A)=-\chi (A)$. This gives rise to
a map
$$G_{[h]}: K\rightarrow \R/\Gamma,$$
where $\Gamma=2\pi\chi(S)\Z$ when $S$ is closed and $\Gamma=0$ when
$S$ has boundary.  When restricted to the Torelli group
$\mathcal{T}(S)$, we have a homomorphism:
$$G:\mathcal{T}(S)\rightarrow H^1(S;\R/\Gamma)\simeq
Hom(H_1(S;\Z),\R/\Gamma),$$
$$[h]\mapsto G_{[h]}.$$
Since the pseudo-Anosov representative of a mapping class is
basically unique, we ask:

\begin{q}
Is $F_{h}=G_{[h]}$ for $h$ pseudo-Anosov and in $\mathcal{T}(S)$, up
to an overall constant factor?
\end{q}

Finally, we briefly discuss the relationship to the {\em
monotonicity condition} for an area-preserving diffeomorphism $h$,
described in Seidel~\cite{Se}.  Suppose that $\chi(S)<0$. On
$\Sigma(S,h)$ consider the tangent bundle $W$ to the fibers and let
$c_1 (W)$ be its first Chern class. The {\em monotonicity condition}
requires that $[\omega_h]=\lambda c_1(W)$ for some real number
$\lambda$. Using the notation from Lemma~\ref{lemma:monotonicity},
one can verify that $\langle c_1(W),C'\rangle=\chi(C)$ for homology
classes of type $C'$. Here $C$ is the surface with $\bdry
C=h(\gamma)-\gamma$. This means that monotonicity holds if and only
if $F_h$ and $G_{[h]}$ are proportional. (A similar, but slightly
more complicated, monotonicity condition also appears in the
definition of {\em periodic Floer homology} of $h$. See~\cite{HS}.)

\s\s\n {\em Acknowledgements.}  We thank Yasha Eliashberg and
Sylvain Gervais for very helpful conversations.  We also thank
Andrew Cotton-Clay and Dan Margalit for their comments on the first
version of the paper.


\begin{thebibliography}{}

\bibitem[Ca]{Ca}
E.\ Calabi, \textit{On the group of automorphisms of a symplectic
manifold},  Problems in analysis (Lectures at the Sympos.\ in honor
of Salomon Bochner, Princeton Univ., Princeton, N.J., 1969), 1--26.
Princeton Univ.\ Press, Princeton, N.J., 1970.

\bibitem[Co]{Co}
V.\ Colin, \textit{Livres ouverts en g\'eom\'etrie de contact
(d'apr\`es Emmanuel Giroux)}, Ast\'erisque \textbf{311}, Expos\'e
969, S\'eminaire Bourbaki, Soci\'et\'e Math\'ematique de France
(2008), 91--117.

\bibitem[CH1]{CH1}
V.\ Colin and K.\ Honda, \textit{Stabilizing the monodromy map of
open books decompositions}, Geom. Dedicata \textbf{132} (2008),
95--103.

\bibitem[CH2]{CH2}
V.\ Colin and K.\ Honda, \textit{Reeb vector fields and open book
decompositions}, preprint 2008.

\bibitem[FLP]{FLP}
A.\ Fathi, F.\ Laudenbach and V.\ Po\'enaru, \textit{Travaux de
Thurston sur les surfaces}, Ast\'erisque \textbf{66-67},
Soci\'et\'e Math\'ematique de France (1991/1971).

\bibitem[HS]{HS}
M.\ Hutchings and M.\ Sullivan, \textit{The periodic Floer homology
of a Dehn twist},  Algebr.\ Geom.\ Topol.\ {\bf 5} (2005), 301--354.

\bibitem[OU]{OU}
J.\ Oxtoby and S.\ Ulam, \textit{Measure-preserving homeomorphisms
and metrical transitivity}, Ann.\ of Math.\ (2) {\bf 42} (1941),
874--920.

\bibitem[Se]{Se}
P.\ Seidel, \textit{Symplectic Floer homology and the mapping class
group}, Pacific J.\ Math.\  {\bf 206}  (2002),  219--229.

\end{thebibliography}
\end{document}